\documentclass[a4paper,reqno]{amsart}
\usepackage{amssymb, amsmath, amscd}
\usepackage{color}

\newtheorem{theorem}{Theorem}[section]
\newtheorem{lemma}[theorem]{Lemma}
\newtheorem{remark}[theorem]{Remark}
\newtheorem{definition}[theorem]{Definition}
\makeatletter
 \@addtoreset{equation}{section}
\makeatother
\begin{document}
\def\KV{homothety}
\title{Homothety Curvature Homogeneity}
\author{E. Garc\'{\i}a-R\'{\i}o, P.  Gilkey  \text{and} S. Nik\v cevi\'c}
\address{EG: Faculty of Mathematics,
University of Santiago de Compostela, Spain}
\email{eduardo.garcia.rio@usc.es}
\address{PG: Mathematics Department, \; University of Oregon, \;\;
  Eugene \; OR 97403 \; USA}
\email{gilkey@uoregon.edu}
\address{SN: Mathematical Institute, Sanu, Knez Mihailova 36, p.p. 367,
11001 Belgrade, Serbia}
\email{stanan@mi.sanu.ac.rs}


\begin{abstract}
We examine the difference between several notions of curvature homogeneity
and show that the notions introduced by Kowalski and Van\v{z}urov\'{a} are
genuine generalizations of the ordinary notion
of $k$-curvature homogeneity.
The homothety group plays an essential role in the analysis.
\end{abstract}
\maketitle

\section{Introduction}


Let $\mathcal{M}=(M,g)$ be a pseudo-Riemannian manifold of
dimension $m\ge3$.
Let $\nabla^kR\in\otimes^{k+4}T^*M$ denote
the $k^{\operatorname{th}}$ covariant derivative of the curvature tensor
and let  $\nabla^k\mathfrak{R}\in
\otimes^{k+2}T^*M\otimes\operatorname{End}(TM)$
denote the $k^{\operatorname{th}}$ covariant derivative
of the curvature operator. These are related by the identity:
$$
\nabla^kR(x_1,x_2,x_3,x_4;x_5,...,x_{k+4})
=g(\nabla^k\mathfrak{R}(x_1,x_2;x_5,...,x_{k+4})x_3,x_4)\,.
$$
\begin{definition}\label{defn-1.1}\rm
$\mathcal{M}$ is said to be {\it k-curvature homogeneous}
if given $P,Q\in M$,  there is a linear isometry
$\phi:T_PM\rightarrow T_QM$ so that
$\phi^*(\nabla^\ell R_Q)=\nabla^\ell R_P$ for $0\le\ell\le k$.
\end{definition}

There is a slightly different version of curvature homogeneity that
we shall discuss here and which,
motivated by the seminal work of
Kowalski and Van\v{z}urov\'{a} \cite{KV11,KV12}, we shall call
{\it \KV\ $k$-curvature homogeneity}.
In Definition~\ref{defn-1.1}, we may replace the curvature tensor $R$
by the curvature operator $\mathfrak{R}$ since we are
dealing with isometries. This is not the case
when we deal with homotheties and the variance is crucial.
We will establish the following result in Section~\ref{sect-2}:
\begin{lemma}\label{lem-1.2}
The following conditions
are equivalent and if any is satisfied, then $\mathcal{M}$ will be said to be
{\bf \KV\ $k$-curvature homogeneous}:
\begin{enumerate}
\item Given any two points $P,Q\in M$, there is a linear homothety
$\phi=\phi_{P,Q}$ from  $T_PM$ to $T_QM$  so that if $0\le\ell\le k$, then
$\phi^*(\nabla^\ell\mathfrak{R}_Q)=\nabla^\ell\mathfrak{R}_P$.
\item Given any two points $P,Q\in M$, there exists a linear isometry
$\Phi=\Phi_{P,Q}$ from $T_PM$ to $T_QM$
and there exists $0\ne\lambda=\lambda_{P,Q}\in\mathbb{R}$
so that if $0\le\ell\le k$, then
$\Phi^*(\nabla^\ell R_Q)=\lambda^{-\ell-2}\nabla^\ell R_P$.
\item There exist constants $\varepsilon_{ij}$ and constants
$c_{i_1...i_{\ell+4}}$ such that
 for all $Q\in M$, there is a basis $\{\xi_1^Q,...,\xi_m^Q\}$ for $T_QM$ and there there exists
 $0\ne\lambda_Q\in\mathbb{R}$ so that if
 $0\le\ell\le k$, then for all $i_1,i_2,...$,
 $g_Q(\xi^Q_{i_1},\xi^Q_{i_2})=\varepsilon_{i_1i_2}$ and
 $$
\nabla^\ell R_Q(\xi^Q_{i_1},\xi^Q_{i_2},\xi^Q_{i_3},\xi^Q_{i_4};\xi^Q_{i_5},...,\xi^Q_{\ell+4})
 =\lambda_Q^{-\ell-2}c_{i_1...i_{\ell+4}}\,.
$$
 \end{enumerate}
 \end{lemma}

\begin{remark}\rm
This agrees with Proposition~0.1 of Kowalski and Van\v{z}urov\'{a} \cite{KV12}.
If we can take $\lambda_{P,Q}=1$ for all $P$ and $Q$, then
$\mathcal{M}$ is $k$-curvature homogeneous.
But we shall see in Theorem~\ref{thm-1.6},
 there are examples which are \KV\ $2$-curvature homogeneous which
 are not $2$-curvature homogeneous and thus $\lambda$ varies with the point.
\end{remark}

Motivated by Lemma~\ref{lem-1.2}, we make the following:
\begin{definition}\rm
A $k$-curvature model is a tuple
$\mathfrak{M}_k:=(V,\langle\cdot,\cdot\rangle,A,...,A^k)$
where $\langle\cdot,\cdot\rangle$
is a non-degenerate inner product on an $m$-dimensional real
vector space $V$ and where
$A^i\in\otimes^{4+i}(V^*)$. We say that two
$k$-curvature models $\mathfrak{M}_k^1$ and $\mathfrak{M}_k^2$ are
{\it \KV\ isomorphic} if there
is a linear isometry $\Phi$ from $(V^1,\langle\cdot,\cdot\rangle^1)$ to
 $(V^2,\langle\cdot,\cdot\rangle^2)$
and if there exists $\lambda\in\mathbb{R}$ so that
$\Phi^*A^{\ell,2}=\lambda^{-\ell-2}A^{\ell,1}$
for $0\le\ell\le k$.
\end{definition}

Lemma~\ref{lem-1.2} shows
a pseudo-Riemannian manifold $\mathcal{M}$ is
\KV\ $k$-curvature homogeneous if and only if
there exists a $k$-curvature model $\mathfrak{M}_k$ so that
$\mathfrak{M}_k$  is \KV\ isomorphic to
$(T_PM,g_P,R_P,...,\nabla^kR_P)$ for all $P$ in $M$.

\subsection{Structure Groups} Let $\mathcal{D}(M)$ denote the group of diffeomorphisms of a pseudo-Riemannian manifold $\mathcal{M}$. We
define the group of isometries $\mathcal{I}(\mathcal{M})$ and the group of homotheties $\mathcal{H}(\mathcal{M})$
by setting:
\begin{eqnarray*}
&&\mathcal{I}(\mathcal{M}):=\{T\in\mathcal{D}(M):T^*g=g\},\\
&&\mathcal{H}(\mathcal{M}):=\{T\in\mathcal{D}(M):
\exists0\ne\lambda=\lambda(T)\in\mathbb{R}:T^*g=\lambda^2g\}\,.
\end{eqnarray*}
We say that $\mathcal{M}$ is {\it homogeneous} if
$\mathcal{I}(\mathcal{M})$ acts transitively on $M$. Similarly,
$\mathcal{M}$ is said to be {\it homothety homogeneous} if
$\mathcal{H}(\mathcal{M})$ acts transitively
on $M$. There are similar local notions where the transformation $T$
is not assumed globally defined.

Homothety homogeneity is essentially a local property. If $(\mathcal{M},g)$ is a complete homothety homogeneous manifold, then there exist $m$-linearly independent homothetic vector fields on $M$.
$(\mathcal{M},g)$ is homogeneous if all of them are Killing and moreover the existence of some non-Killing homothetic vector fields is very restrictive. A complete Riemannian manifold which admits a non-Killing homothetic vector field must be flat \cite{Tashiro}, and hence it follows that a non-flat complete homothety homogeneous manifold is necessarily homogeneous in the Riemannian setting. The situation is not so rigid in the Lorentzian case where pp-wave metrics support non-Killing
homothetic vector fields (see for example \cite{Alek, Ku-Ra, Steller} and references therein).

\subsection{Stability}
Assertion~(1) in
the following result was established by Singer \cite{S60} in the
Riemannian context and by Podesta and Spiro \cite{PS96}
in the pseudo-Riemannian setting. In Section~\ref{sect-3}, we
will use results of \cite{PS96} to establish Assertion~(2) which
extends these results to the \KV~ setting. Recall that the linear orthogonal
group $\mathcal{O}$ and the linear homothety group $\mathcal{H}$
in dimension $m$ but arbitrary signature satisfy:
$$
\dim\{\mathcal{O}\}:=\textstyle\frac12m(m-1)\text{ and }
\dim\{\mathcal{H}\}=\textstyle\frac12m(m-1)+1\,.$$

\begin{theorem}\label{thm-1.5}
Let $\mathcal{M}=(M,g)$ be a pseudo-Riemannian manifold.
\begin{enumerate}
\item The following Assertions are equivalent:
\begin{enumerate}
\item $\mathcal{M}$ is locally homogeneous.
\item $\mathcal{M}$ is $k$-curvature homogeneous for all $k$.
\item $\mathcal{M}$ is $k$-curvature homogeneous for $k=\frac12m(m-1)$.
\end{enumerate}
\item The following Assertions are equivalent:
\begin{enumerate}
\item $\mathcal{M}$ is locally homothety homogeneous.
\item $\mathcal{M}$ is \KV\ $k$-curvature homogeneous for all $k$.
\item $\mathcal{M}$ is \KV\ $k$-curvature homogeneous for
$k=\frac12m(m-1)+1$.
\end{enumerate}
\end{enumerate}
\end{theorem}

\subsection{Homothety homogeneous manifolds
that are not $0$-curvature homogeneous}
Let $m\ge3$ and let
$\mathcal{N}=(N,g_N)$ be a homogeneous
pseudo-Riemannian manifold of dimension $m-1$. Set:
$$\mathcal{M}_t:=
(\mathbb{R}\times N,g_{M,t})\text{ where }g_{M,t}:=e^{tx}(dx^2+g_N)\,.$$
Let $\tau_{\mathcal{N}}$ and $\tau_{\mathcal{M}_t}$ denote the
scalar curvature of $\mathcal{N}$ and of $\mathcal{M}_t$, respectively.
We will establish the following result in Section~\ref{sect-4}:

\begin{theorem}\label{thm-1.6}
Assume that $m\ge3$.
\begin{enumerate}
\item $\mathcal{M}_t$ is homothety homogeneous and hence \KV\ $k$-curvature homogeneous
for all $k$.
\item Suppose that $t\ne0$ is chosen so that
$\{\tau_{\mathcal{N}}-\frac{(m-1)(m-2)}4t^2\}\ne0$.
Then $\mathcal{M}_t$ is not $0$-curvature homogeneous and in particular
not locally homogeneous.
\end{enumerate}
\end{theorem}

Let $\lambda(\Phi)$ be the homothety constant so that
$\Phi^*(g)=\lambda^2(\Phi)g$; we may always assume that
$\lambda(\Phi)>0$. The manifolds of Theorem~\ref{thm-1.6} are
cohomogeneity one, i.e. the group of isometries acts
transitively on a family of hyper surfaces which foliate the manifold. This
is in fact the general setting as we shall show in Section~\ref{sect-5}:
\begin{theorem}\label{thm-1.7}
Assume $\mathcal{M}=(M,g)$ is a
pseudo-Riemannian
manifold which is not homogeneous, which is homothety homogeneous,
and which has $|R|^2\ne0$. Fix a base point $P_0$ of $M$ and
define a smooth function $\mu\in C^\infty(M)$ by setting:
$$\mu(P):=\frac{|R|^2_{P_0}}{|R|^2_P}\,.$$
\begin{enumerate}
\item If $\Phi_1$ and $\Phi_2$ are homotheties, then
$\lambda^2(\Phi_1\circ\Phi_2)=\lambda^2(\Phi_1)\cdot\lambda^2(\Phi_2)$.
\item $\lambda^2(\Phi)=\frac{\mu(\Phi(P))}{\mu(P)}$
for any $P\in M$.
\item $d\mu\ne0$.
\item Let $M_c:=\{P\in M:\mu(P)=c\}$
define smooth submanifolds of $M$. Given any points $P_i\in M_c$,
there exists an isometry of $M$ which preserves $M_c$ so that
$\Phi(P_1)=\Phi(P_2)$. Thus $(M,g)$ has
cohomogeneity one.
\end{enumerate}
\end{theorem}

Our analysis is local; if $\mathcal{M}$ is only assumed to be locally
homothety homogeneous, then we may conclude $\mathcal{M}$ is
locally cohomogeneity one.
In the Riemannian setting, since $\mathcal{M}$ is not homogeneous,
it is not flat and hence the condition $|R|^2\ne0$ is automatic. 
In the
higher signature setting, there are manifolds which are not flat but which
satisfy $|R|^2=0$ and, more generally, have all their Weyl scalar invariants
vanish. These are called VSI manifolds
-- we refer to \cite{AMCH12,CMPP04} for further
details. There is a vast literature concerning VSI manifolds.

\subsection{Walker Lorentzian 3 dimensional manifolds}\label{sect-1.4}
Section~\ref{sect-6} is devoted to the study of a very specific
family of examples.
Let $\mathcal{M}=(M,g_M)$ be a $3$-dimensional Lorentzian manifold
which admits a parallel null
vector field, i.e. $\mathcal{M}$ is a $3$-dimensional Walker manifold.
Such a manifold admits local
adapted coordinates $(x,y,\tilde x)$ so that the (possibly) non-zero components of the metric are given by
$$g(\partial_x,\partial_x)=-2f(x,y),\quad g(\partial_x,\partial_{\tilde x})=g(\partial_y,\partial_y)=1\,.$$

We shall denote this manifold by $\mathcal{M}_f$. We have (see \cite{GNS12} Theorem 2.10 and Theorem 2.12):

\begin{theorem}\label{thm-1.8}
Suppose that $f_{yy}$ is never zero.
\begin{enumerate}
\item $\mathcal{M}_f$ is 1-curvature homogeneous if and only if exactly
one of the following three possibilities holds:
\begin{enumerate}
\item $f_{yy}(x,y)=ay^2$ for $a\ne0$. This manifold is symmetric.
\item $f_{yy}(x,y)=\alpha(x)e^{by}$ where $0\ne b\in\mathbb{R}$ and where $\alpha(x)$ is arbitrary.
\item $f_{yy}(x,y)=\alpha(x)$ where $\alpha(x)=c\cdot\alpha_x^{3/2}$ for some $c\ne0\in\mathbb{R}$.
\end{enumerate}
\item $\mathcal{M}$ is $2$-curvature homogeneous if and only if it falls into one of the three families, all of which are locally homogeneous:
\begin{enumerate}
\item $f=b^{-2}\alpha(x)e^{by}+\eta(x)y+\gamma(x)$ where
$0\ne b\in\mathbb{R}$, where
 $\alpha(x)\ne0$, and where
$\eta(x)=b^{-1}\alpha^{-1}(x)
\{\alpha_{xx}(x)-\alpha_x^2(x)\alpha(x)^{-1}\}$.
\item $f=a(x-x_0)^{-2}y^2+\beta(x)y+\gamma(x)$ where
$0\ne a\in\mathbb{R}$.
\item $f=ay^2+\beta(x)y+\gamma(x)$ where $0\ne a\in\mathbb{R}$.
\end{enumerate}
\end{enumerate}\end{theorem}

We will establish the following analogue of Theorem~\ref{thm-1.8}
for \KV\ curvature homogeneity in Section~\ref{sect-6}:

\begin{theorem}\label{thm-1.9}
Suppose $f_{yy}$ is never zero and non-constant.\begin{enumerate}
\item If $f_{yyy}$ never vanishes, then
$\mathcal{M}_f$ is \KV\ 1-curvature
homogeneous.
\item If $f_{yy}=\alpha(x)$ with $\alpha_x$ never zero, then $\mathcal{M}_f$
is \KV\ 1-curvature homogeneous if and only if
$f=a(x-x_0)^{-2}y^2+\beta(x)y+\gamma(x)$
where $0\ne a\in\mathbb{R}$. This manifold is locally homogeneous.
\item Assume that $\mathcal{M}_f$ is \KV\ 2-curvature homogeneous,
and that $f_{yy}$ and $f_{yyy}$ never vanish.
Then $\mathcal{M}_f$ is locally isometric to one of the examples:
\begin{enumerate}
\item $\mathcal{M}_{\pm e^{ay}}$ for some $a\ne0$ and
$M=\mathbb{R}^3$. This manifold is
homogeneous.
\item $\mathcal{M}_{\pm\ln(y)}$ for some $a\ne0$ and $M=\mathbb{R}\times(0,\infty)\times\mathbb{R}$. This manifold is homothety homogeneous
but not locally homogeneous.
\item $\mathcal{M}_{\pm y^\varepsilon}$ for $\varepsilon\ne0,1,2$ and
$M=\mathbb{R}\times(0,\infty)\times\mathbb{R}$. This manifold is homothety homogeneous
but not locally homogeneous.
\end{enumerate}\end{enumerate}
\end{theorem}

It will follow from our analysis that the manifolds
$\mathcal{M}_{\pm\ln(y)}$ and $\mathcal{M}_{\pm y^c}$
are homothety homogeneous VSI manifolds which are cohomogeneity one,
thereby exhibiting non-trivial examples in the VSI setting.
We also refer to recent work of Dunn and McDonald \cite{DuMc2013}
for related work on homothety curvature homogeneous manifolds.

\subsection{Variable \KV\ curvature homogeneity}
In fact, the definition we have used in this paper differs subtly but in an
important fashion from that originally given by Kowalski and Van\v{z}urov\'{a} \cite{KV12};
in that paper the scaling constant $\lambda$ was permitted to depend
on $\ell$ and this gives rise to the notion of
{\it variable \KV\ $k$-curvature homogeneity}. There are 4 different
definitions which may be summarized as follows; we repeat two of
the definitions to put the new definitions in context:
\begin{definition}\label{defn-1.10}
\rm
Let $\mathfrak{M}_k:=\{V,\langle\cdot,\cdot\rangle,A_0,...,A_k\}$ be
a $k$-curvature model and let $\mathcal{M}=(M,g)$ be a pseudo-Riemannian manifold.
\begin{enumerate}
\item Static isometries that are independent of $k$. Recall that:
\begin{enumerate}
\item $\mathcal{M}$ is {\it $k$-curvature homogeneous with model $\mathfrak{M}_k$}
 if for any $P$ in $M$, there exists an isometry
 $\phi_P:T_PM\rightarrow V$ so
$\phi_P^*A_\ell=\nabla^\ell R_P$ for $0\le\ell\le k$.
\item $\mathcal{M}$ is {\it\KV\ $k$-curvature homogeneous with model $\mathfrak{M}_k$}
 if for any $P$ in $M$, there is an isometry $\phi_P:T_PM\rightarrow V$
 and a scaling factor $0\ne\lambda\in\mathbb{R}$ so that
$\phi_P^*A_\ell=\lambda^{2+\ell}R_P$ for $0\le\ell\le k$.
\end{enumerate}
\item Variable isometries that depend on $k$. We shall say that:
\begin{enumerate}
\item $\mathcal{M}$ is
{\it variable-$k$-curvature homogeneous with model $\mathfrak{M}_k$}
if for every $P$ in $M$ and if for every $0\le\ell\le k$,
there exist isometries $\phi_{P,\ell}:T_PM\rightarrow V$ so that
$\phi_{P,\ell}^*A_\ell=\nabla^\ell R_P$.
\item We say that $\mathcal{M}$ is
{\it variable \KV\ $k$-curvature homogeneous with model
$\mathfrak{M}_k$} if for every$P$ in $M$ and if for every $0\le\ell\le k$,
there are isometries $\phi_{P,\ell}:T_PM\rightarrow V$ and scaling factors
$0\ne\lambda_\ell\in\mathbb{R}$ so that
$\phi_{P,\ell}^*A_\ell=\lambda_\ell^{2+\ell}R_P$.
\end{enumerate}
\end{enumerate}\end{definition}
In Section~\ref{sect-7}, we will use the examples which were studied
in Section~\ref{sect-6} to show that Theorem~\ref{thm-1.5} fails in the
context of variable curvature homogeneity and hence also for variable
\KV\ curvature homogeneity:
\begin{theorem}\label{thm-1.11}
\ \begin{enumerate}
\item Let $k$ be a positive integer.
There exists $f_k$ so that $\mathcal{M}_{f_k}$ is
variable $k$-curvature homogeneous but not variable $k+1$-curvature homogeneous.
\item The manifold $\mathcal{M}_{\frac12e^xy^2}$ is
variable $k$-curvature homogeneous for all $k$, but is not \KV\ 1-curvature homogeneous and hence not locally homogeneous.
\item In Definition~\ref{defn-1.10} we have the following implications:
$$
(1a)\quad\Rightarrow\quad(1b)\quad\Rightarrow\quad(2b)\quad\text{and}\quad
(1a)\quad\Rightarrow\quad(2a)\quad\Rightarrow\quad(2b)\,.
$$
All other possible implications are false.
\end{enumerate}
\end{theorem}

\section{The proof of Lemma~\ref{lem-1.2}}\label{sect-2}

Assume that Assertion~(1) of Lemma~\ref{lem-1.2} holds.
This means that given any two points $P$ and $Q$ in $M$,
there exists a linear homothety
$\phi:T_PM\rightarrow T_QM$ so that if $0\le\ell\le k$
and if $\{x_i\}$ are vectors in $T_PM$, then we have that:
\begin{eqnarray*}
&&g_Q(\phi x_1,\phi x_2)=\lambda^2g_P(x_1,x_2),\text{ and}\\
&&\phi\left\{\nabla^\ell\mathfrak{R}_P(x_1,x_2;x_5,x_{\ell+4})x_3\right\}
=\nabla^\ell\mathfrak{R}_Q(\phi x_1,\phi x_2;\phi x_5,...,\phi x_{\ell+4})\phi x_3\,.
\end{eqnarray*}
Taking the inner product with $x_4$ permits us to rewrite the
second condition, which involves the curvature operator,
in terms of the curvature tensor:
\begin{eqnarray*}
&&\lambda^2\nabla^\ell R_P(x_1,x_2,x_3,x_4;x_5,...,x_{\ell+4})\\
&=&\lambda^2g_P\left(\nabla^\ell\mathfrak{R}_P(x_1,x_2;x_5...,x_{\ell+4})x_3,x_4\right)\\
&=&g_Q(\phi\nabla^\ell\mathfrak{R}_P(x_1,x_2;x_5...,x_{\ell+4})x_3,\phi x_4)\\
&=&g_Q(\nabla^\ell\mathfrak{R}_Q(\phi x_1,\phi x_2;\phi x_5,...,\phi x_{\ell+4})\phi x_3,\phi x_4)\\
&=&\nabla^\ell R_Q(\phi x_1,\phi x_2,\phi x_3,\phi x_4;\phi x_5,...,\phi x_{\ell+4})\,.
\end{eqnarray*}
We set $\Phi:=\lambda^{-1}\phi$. We can rewrite these equations in the form:
\begin{eqnarray*}
&&g_Q(\Phi x_1,\Phi x_2)=\lambda^{-2}g_Q(\phi x_1,\phi x_2)=g_P(x_1,x_2),\\
&&\lambda^2\nabla^\ell R_P(x_1,x_2,x_3,x_4;x_5,...,x_{\ell+4})\\
&=&\nabla^\ell R_Q(\phi x_1,\phi x_2,\phi x_3,\phi x_4;\phi x_5,...,\phi x_{\ell+4})\\
&=&\lambda^{\ell+4}\nabla^\ell R_Q(\Phi x_1,\Phi x_2,\Phi x_3,\Phi x_4;\Phi x_5,...,\Phi x_{\ell+4})\\
&=&\lambda^{\ell+4}\Phi^*(\nabla^\ell R_Q)(x_1,x_2,x_3,x_4;x_5,...,x_{\ell+4})\,.
\end{eqnarray*}
This shows $\Phi$ is an isometry from $T_PM$ to $T_QM$ so
$\Phi^*(\nabla^\ell R_Q)=\lambda^{-2-\ell}\nabla^\ell R_P$. Consequently in
Lemma~\ref{lem-1.2}, Assertion~(1) $\Rightarrow$
Assertion~(2);
the proof of the converse implication is similar and will be omitted.

Suppose that Assertion~(2) of Lemma~\ref{lem-1.2} holds.
Fix a base point $P\in M$ and fix a basis $\{\xi_1^P,...,\xi_m^P\}$
for $T_PM$. Set $\varepsilon_{ij}:=g_P(\xi_i,\xi_j)$ and
$c_{i_1,...,i_{\ell+4}}:=\nabla^\ell R(\xi_{i_1},...,\xi_{i_{\ell+4}})$.
Let $Q\in M$. By assumption,
there is an isometry $\Phi$ from $T_PM$ to $T_QM$ so:
$$\Phi^*(\nabla^\ell R_Q)=\lambda_{Q}^{-\ell-2}\nabla^\ell R_P\text{ for }0\le\ell\le k\,.$$
Set $\xi_i^Q:=\Phi\xi_i^P$. Then:
\begin{eqnarray*}
&&g_Q(\xi_i^Q,\xi_j^Q)=g_Q(\Phi\xi_i^P,\Phi\xi_j^P)=g_P(\xi_i^P,\xi_j^P)=\varepsilon_{ij},\\
&&\nabla^\ell R_Q(\xi_{i_1}^Q,...,\xi_{i_{\ell+4}}^Q)=
\nabla^\ell R_Q(\Phi\xi_{i_1}^P,...,\Phi\xi_{i_{\ell+4}}^P)\\
&&\qquad=
\lambda_{Q}^{-\ell-2}\nabla^\ell R_P(\xi_{i_1}^P,...,\xi_{i_{\ell+4}}^P)=
\lambda_{Q}^{-\ell-2}c_{i_1,...,i_{\ell+4}}\,.
\end{eqnarray*}
This shows in Lemma~\ref{lem-1.2} that Assertion~(2)
$\Rightarrow$ Assertion~(3); the proof of the converse implication is
similar and will be omitted.
\hfill\qed

\section{The proof of Theorem~\ref{thm-1.5}}\label{sect-3}

Fix the signature $(p,q)$ throughout this section where $m=p+q$.
Let $(V,\langle\cdot,\cdot\rangle)$ be an inner product space
of signature $(p,q)$.
Let $\mathcal{O}$ and $\mathcal{HO}$
be the orthogonal group and homothety group of signature $(p,q)$
respectively;
\begin{eqnarray*}
&&\mathcal{O}:=\{T\in\operatorname{GL}(V):
T^*\langle\cdot,\cdot\rangle=\langle\cdot,\cdot\rangle\},\\
&&\mathcal{HO}:=\{T\in\operatorname{GL}(V):
T^*\langle\cdot,\cdot\rangle=\lambda^2\langle\cdot,\cdot\rangle\}
\text{ for some }0<\lambda\in\mathbb{R}\}\,.
\end{eqnarray*}
Note that $\mathcal{HO}=\mathbb{R}^+\times\mathcal{O}$.
Fix a basis $\{v_1,...,v_m\}$ for $V$
and let $\epsilon_{ij}:=\langle v_i,v_j\rangle$. Let $\mathcal{M}=(M,g)$
be a pseudo-Riemannian manifold of signature $(p,q)$. Let
$\mathcal{F}(M)$ be the bundle of frames $\vec u=(u_1,...,u_m)$ for
$TM$.
Let $\mathcal{O}(\mathcal{M})$ and $\mathcal{HO}(\mathcal{M})$
be the sub-bundles:
\begin{eqnarray*}
&&\mathcal{O}(\mathcal{M}):=\{\vec u\in\mathcal{F}(M):
g(u_i,u_j)=\epsilon_{ij}\},\\
&&\mathcal{HO}(\mathcal{M}):=\{\vec u\in\mathcal{F}(M):
g(u_i,u_j)=\lambda^\epsilon_{ij}
\text{ for some }0<\lambda\in\mathbb{R}\}\,.
\end{eqnarray*}
Then $\mathcal{O}(M)$ is a principal $\mathcal{O}$ bundle while
$\mathcal{HO}(M)$ is a principal $\mathcal{HO}$ bundle.
We use the Levi-Civita connection, which is invariant under homothetic
transformations, to define a $\mathcal{HO}$ structure on $\mathcal{M}$
with a canonical connection.

Let $\mathfrak{ho}(T_P\mathcal{M})$ be the Lie algebra of the group of
homothetic transformations of $T_P\mathcal{M}$. Let
$$
s_0:=\dim\{\mathfrak{ho}(T_P\mathcal{M})\}=\frac12m(m-1)+1\,.
$$
For $0\le s\le s_0$, we consider the subalgebras defined by:
\begin{eqnarray*}
&&\mathfrak{ho}^0(T_P\mathcal{M})
=\{ a\in\mathfrak{ho}(T_P\mathcal{M})\,;\,\, a\cdot R=0\},\\
&&\mathfrak{ho}^s(T_P\mathcal{M})
=\{ a\in\mathfrak{ho}^{s-1}(T_P\mathcal{M})\,;\,\, a\cdot \nabla^sR=0\}\,.
\end{eqnarray*}
Clearly
$$
\mathfrak{ho}^s(T_P\mathcal{M})\supset
\mathfrak{ho}^{s+1}(T_P\mathcal{M})\text{ for all }s$$
so we have a decreasing sequence of subalgebras of
$\mathfrak{ho}(T_P\mathcal{M})$.
Let the {\it Singer number} $s(P)$ be the first integer stabilizing this
sequence above, i.e.:
$$
\mathfrak{ho}^{s(P)+r}(T_P\mathcal{M})
=\mathfrak{ho}^{s(P)}(T_P\mathcal{M})\text{ for all }r\geq 1\,.
$$
Now the assumption that $(\mathcal{M},g)$ is
\KV~$k$-curvature homogeneous for some $k\ge\frac12m(m-1)+1$
shows that the Singer number $s(P)$ is constant on $M$.
The equivalences of Theorem~\ref{thm-1.5}~(2) now
follow from the work of Podesta and Spiro \cite{PS96}.
(See also \cite{Op} for an extension to the affine setting.)

\section{The proof of Theorem~\ref{thm-1.6}}\label{sect-4}

Let $\mathcal{N}=(N,g_N)$ be a homogeneous pseudo-Riemannian manifold
of dimension $m-1\ge2$
and let $\mathcal{M}=(\mathbb{R}^\times N,g_{M,t})$ where we take
$g_{M,t}:=e^{tx}(dx^2+g_N)$.
Let $T_a(x):=x+a$ and let $\theta$ be
 an isometry of $\mathcal{N}$. Then:
 $$
 (T_a\times\theta)^*(g_{M,t})=T_a^*(f_t)(T_a\times\theta)^*(g_{M,0})
=e^{ta}g_{M,t}\,.
 $$
This shows that $\mathcal{T}\subset\mathcal{H}(\mathcal{M}_t)$.
Elements of this form act transitively on $M$ and hence
$\mathcal{M}_t$ is  homothety homogeneous.

We now examine the curvature tensor. Fix $t$ and fix a point $P\in N$.
Let $g=g_N$ and let $\tilde g=g_{M,t}$.
Choose local coordinates
$y=(y^1,...,y^{m-1})$ centered at $P$.
Let indices $u,v,w$ range from $0$ to $m-1$
and index the coordinate frame
$(\partial_x,\partial_{y_1},...,\partial_{y_{m-1}})$;
indices $i,j,k$ range from $1$ to $m-1$ and index the
coordinate frame $(\partial_{y_1},...,\partial_{y_{m-1}})$.
Let $\Gamma$ be the Christoffel symbols of $g$ and $\tilde\Gamma$
be the Christoffel symbols of $\tilde g$. Let $\delta_i^j$ be the Kronecker
index. We compute:
$$\begin{array}{llll}
\tilde g_{00}=e^{tx},&\tilde g_{0i}=0,&\tilde g_{ij}=e^{tx}g_{ij},\\
\tilde \Gamma_{000}=\frac12te^{tx}&
\tilde \Gamma_{0i0}=0,&
\tilde \Gamma_{ij0}=-\frac12te^{tx}g_{ij},
\vphantom{\vrule height 12pt}\\
\tilde \Gamma_{00i}=0,&
\tilde\Gamma_{0ij}=\frac12te^{tx}g_{ij},&
\tilde \Gamma_{ijk}=e^{tx}\Gamma_{ijk},
\phantom{\vrule height 12pt}\\
\tilde \Gamma_{00}{}^0=\frac12t&
\tilde \Gamma_{0i}{}^0=0,&
\tilde \Gamma_{ij}{}^0=-\frac12tg_{ij},
\vphantom{\vrule height 12pt}\\
\tilde \Gamma_{00}{}^i=0,&
\tilde\Gamma_{0i}{}^j=\frac12t\delta_i^j,&
\tilde \Gamma_{ij}{}^k=\Gamma_{ij}{}^k.
\phantom{\vrule height 12pt}
\end{array}$$
Thus the covariant derivatives are given by
$$\begin{array}{lll}
\tilde\nabla_{\partial_x}\partial_x=\frac12t\partial_x,&
\tilde\nabla_{\partial_x}\partial_{y_i}=\frac12t\partial_{y_i},\\
\tilde\nabla_{\partial_{y_i}}\partial_x=\frac12t\partial_{y_i},&
\tilde\nabla_{\partial_{y_i}}\partial_{y_j}
=\Gamma_{ij}{}^k\partial_{y_k}-\frac12tg_{ij}\partial_x.
\vphantom{\vrule height 12pt}
\end{array}$$
We choose the coordinate system so the first derivatives of
$g_{ij}$ vanish at $P$ and hence $\Gamma(P)=0$.
Consequently the curvature operator at $P$ is given by:
\medbreak\quad
$\tilde{\mathcal{R}}(\partial_x,\partial_{y_i})\partial_x=
\{\tilde\nabla_{\partial_x}\tilde\nabla_{\partial_{y_i}}
-\tilde\nabla_{\partial_{y_i}}\tilde\nabla_{\partial_x}\}\partial_x
=\frac12t\tilde\nabla_{\partial_x}\partial_{y_i}
-\frac12t\tilde\nabla_{\partial_{y_i}}\partial_x=0$,
\medbreak\quad
$\tilde{\mathcal{R}}(\partial_x,\partial_{y_i})\partial_{y_j}=
\{\tilde\nabla_{\partial_x}\tilde\nabla_{\partial_{y_i}}
-\tilde\nabla_{\partial_{y_i}}\tilde\nabla_{\partial_x}\}\partial_{y_j}$
\smallbreak\qquad\qquad\qquad\quad
$=\tilde\nabla_{\partial_x}
\{\Gamma_{ij}{}^k\partial_{y_k}-\frac12tg_{ij}\partial_x\}
-\frac12t\tilde\nabla_{\partial_{y_i}}\partial_{y_j}$
\smallbreak\qquad\qquad\qquad\quad
$=\frac12t\{\Gamma_{ij}{}^k\partial_{y_k}-\frac12tg_{ij}\partial_x\}
-\frac12t\{\Gamma_{ij}{}^k\partial_{y_k}-\frac12tg_{ij}\partial_x\}=0$,
\medbreak\quad
$\tilde{\mathcal{R}}(\partial_{y_i},\partial_{y_j})\partial_x=
\{\tilde\nabla_{\partial_{y_i}}\tilde\nabla_{\partial_{y_j}}
-\tilde\nabla_{\partial_{y_j}}\tilde\nabla_{\partial_{y_i}}\}\partial_x
=\frac12t\tilde\nabla_{\partial_{y_i}}\partial_{y_j}-
\frac12t\tilde\nabla_{\partial_{y_j}}\partial_{y_i}=0$,
\medbreak\quad
$\tilde{\mathcal{R}}(\partial_{y_i},\partial_{y_j})\partial_{y_k}=
\{\tilde\nabla_{\partial_{y_i}}\tilde\nabla_{\partial_{y_j}}
-\tilde\nabla_{\partial_{y_j}}\tilde\nabla_{\partial_{y_i}}\}\partial_{y_k}$
\smallbreak\qquad\qquad\qquad\quad
$=\tilde\nabla_{\partial_{y_i}}(\Gamma_{jk}{}^\ell\partial_{y_\ell}
-\frac12tg_{jk}\partial_x)
-\tilde\nabla_{\partial_{y_j}}(\Gamma_{ik}{}^\ell\partial_{y_\ell}
-\frac12tg_{ik}\partial_x)$
\smallbreak\qquad\qquad\qquad\quad
$=R_{ijk}{}^\ell\partial_{y_\ell}-\frac14t^2g_{jk}\partial_{y_i}
+\frac14t^2g_{ik}\partial_{y_j}$.
\medbreak\noindent We can now express the scalar curvature and Ricci
tensor $\{\tilde\rho,\tilde\tau\}$ for $\tilde g$ in terms of the scalar
curvature and Ricci tensor $\{\rho,\tau\}$ for $g$:
$$
\tilde\rho=\rho-\textstyle\frac{m-2}4t^2g\text{ and }
\tilde\tau=e^{-tx}\{\tau-\frac{(m-1)(m-2)}4t^2\}\,.
$$
Suppose that $\tau-\frac{(m-1)(m-2)}4t^2\ne0$
and $t\ne0$. It then follows that $\tilde\tau$ is not constant and
hence $\mathcal{M}_t$ is not $0$-curvature homogeneous. This
completes the proof of Theorem~\ref{thm-1.6}.\hfill\qed

\section{The proof of Theorem~\ref{thm-1.7}}\label{sect-5}
We establish Theorem~\ref{thm-1.7}~(1) by computing:
\begin{eqnarray*}
&&g_{\Phi_1(\Phi_2(P))}((\Phi_1\Phi_2)_*X,(\Phi_1\Phi_2)_*Y)
=\lambda^2(\Phi_1\circ\Phi_2)g_{P}(X,Y)\\
&=&\lambda^2(\Phi_1)g_{\Phi_2(P)}((\Phi_2)_*X,(\Phi_2)_*Y)
=\lambda^2(\Phi_1)\lambda^2(\Phi_2)g_P(X,Y)\,.
\end{eqnarray*}

Let $\{e_i\}$ be a local orthonormal basis for $T_PM$. Then
$\{E_i:=\lambda(\Phi)^{-1}\Phi_*e_i\}$ is an
orthonormal basis for $T_{\Phi(P)}$M.
Since $\Phi$ is a homothety, it commutes with the Levi-Civita connection
and with the curvature operator so we have:
\begin{eqnarray*}
&&R_{\Phi(P)}(E_i,E_j,E_k,E_l)
=\lambda(\Phi)^{-4}g_{\Phi(P)}(\mathcal{R}(\Phi_*e_i,\Phi_*e_j)
\Phi_*e_i,\Phi_*e_l)\\
&=&\lambda(\Phi)^{-4}g_{\Phi(P)}
(\Phi_*\mathcal{R}(e_i,e_j)e_k,\Phi_*e_l)
=\lambda(\Phi)^{-2}g_P(\mathcal{R}(e_i,e_j)e_k,e_l)\\
&=&\lambda(\Phi)^{-2}R_P(e_i,e_j,e_k,e_l)\,.
\end{eqnarray*}
This implies
$$\lambda^2(\Phi)=\frac{|R_P|^2}{|R_{\Phi(P)}|^2}=
\frac{|R|_{P_0}^2}{|R|_{\Phi(P)}^2}\cdot\frac{|R|_P^2}{|R_{P_0}^2}
=\frac{\mu(\Phi(P))}{\mu(P)}\text { for  any }P\in M\,.
$$
Taking $P=P_0$ then yields $\lambda^2(\Phi)=\mu(\Phi(P_0))$
since $\mu(P_0)=1$. Choose $\Phi_1$ so $\Phi_1P_0=P$.
We have:
\begin{eqnarray*}
&&\lambda^2(\Phi)\mu(P)=\lambda^2(\Phi)\mu(\Phi_1P_0)=
\lambda^2(\Phi)\lambda^2(\Phi_1)\\
&=&\lambda^2(\Phi\Phi_1)
=\mu(\Phi\Phi_1P_0)=\mu(\Phi(P))\,.
\end{eqnarray*}
Since $P$ was arbitrary, this shows that we have the intertwining formula
\begin{equation}\label{eqn-5.a}
\Phi^*(\mu)=\lambda(\Phi)\cdot\mu\,.
\end{equation}
Thus $d\mu(\Phi(P))=0$ implies $d\mu(P)=0$.
If $d\mu$ vanishes everywhere, then $\mu$ is constant.
Since $\mu(P_0)=1$, this implies $\mu\equiv1$ so
every homothety is an isometry and $(M,g)$ is homogeneous,
contrary to our assumption. Thus there exists some point where $d\mu\ne0$.
Since the homotheties
act transitively on $M$, Equation~(\ref{eqn-5.a}) implies
 $d\mu$ never vanishes. Consequently,
 the level sets are smooth submanifolds of
$M$. Let $P_i\in M_c$. Choose a  homothety $\Phi$
so $\Phi(P_1)=P_2$. Then
$\lambda^2(\Phi)=1$ so $\Phi$ is an isometry. That means
$\mu(\Phi(P))=\mu(P)$ for all $P$ and hence $\Phi$ preserves all the level sets $M_c$. Since $\Phi$ is an
ambient isometry, $\Phi$ restricts to an isometry of each level set.
This completes the proof of Theorem~\ref{thm-1.7}.\hfill\qed

\medskip
The cohomogeneity one in Theorem~\ref{thm-1.7} also follows in the Riemannian setting from the discussion by Console and Olmos \cite{Co-Ol}, since it follows from Lemma \ref{lem-1.2} and Assertion (2) in Theorem~\ref{thm-1.7} that the regular level sets of the Weyl scalar invariants define a foliation of codimension one unless the manifold is flat.

\section{The proof of Theorem~\ref{thm-1.9}}\label{sect-6}

As the theory is local, we let $\mathcal{M}_f=(\mathbb{R}^3,g_f)$
where
$$
g_f(\partial_x,\partial_x)=-2f(x,y),\quad
g_f(\partial_x,\partial_{\tilde x})=g_f(\partial_y,\partial_y)=1\,.
$$
We suppress the subscript ``$f$" when no confusion is likely to ensure.
We follow the discussion in \cite{GNS12}.
The (possibly) non-zero covariant derivatives are given by:
$$\nabla_{\partial_x}\partial_x=-f_x\partial_{\tilde x}
+f_y\partial_y\text{ and }
\nabla_{\partial_x}\partial_y=\nabla_{\partial_y}\partial_x
=-f_y\partial_{\tilde x}\,.$$
The (possibly) non-zero curvatures and covariant derivatives to order $2$ are:
\begin{eqnarray*}
&&R(\partial_x,\partial_y,\partial_y,\partial_x)=f_{yy},\\
&&\nabla R(\partial_x,\partial_y,\partial_y,\partial_x;\partial_x)=f_{xyy},\\
&&\nabla R(\partial_x,\partial_y,\partial_y,\partial_x;\partial_y)=f_{yyy},\\
&&\nabla R(\partial_x,\partial_y,\partial_y,\partial_x;\partial_x,
\partial_x)=f_{xxyy}-f_yf_{yyy},\\
&&\nabla^2R(\partial_x,\partial_y,\partial_y,\partial_x;\partial_x,\partial_y)=
\nabla^2R(\partial_x,\partial_y,\partial_y,\partial_x;\partial_y,
\partial_x)=f_{xyyy},\\
&&\nabla^2R(\partial_x,\partial_y,\partial_y,\partial_x;\partial_y,
\partial_y)=f_{yyyy}.
\end{eqnarray*}

If $f_{yy}$ vanishes identically, then $\mathcal{M}$ is flat.
The vanishing of $f_{yy}$ is an invariant
of the \KV\ 0-model. Since we are interested in \KV\ curvature homogeneity,
we shall assume $f_{yy}$
never vanishes; since $M$ is connected, either $f_{yy}$ is always positive or
$f_{yy}$ is always negative.
We shall usually assume $f_{yy}>0$ as the other case is handled similarly.
The simultaneous vanishing of $f_{yyy}$ and of $f_{xyy}$ is an invariant of
the \KV\ 1-model. The case $f_{yy}=ay^2$ for $0\ne a\in\mathbb{R}$
gives rise to a symmetric space. We shall therefore assume $f_{yy}$
non-constant. This gives rise to
two cases $f_{yyy}$ never zero and $f_{yyy}$ vanishing identically
but $f_{xyy}$ never zero.

\begin{definition}
\rm Let $\mathfrak{M}_{1,c_1,c_2}$ be the 1-curvature model whose (possibly) non-zero components
are defined by $\varepsilon_{12}=\varepsilon_{22}=1$, $c_{1221}=1$, $c_{12211}=c_1$, and
$c_{12212}=c_2$:
$$\begin{array}{lll}
\langle\xi_1,\xi_2\rangle=1,&\langle\xi_2,\xi_2\rangle=1,&
R(\xi_1,\xi_2,\xi_2,\xi_1)=1,\\
\nabla R(\xi_1,\xi_2,\xi_2,\xi_1;\xi_1)=c_1,&
\nabla R(\xi_1,\xi_2,\xi_2,\xi_1;\xi_2)=c_2.
\end{array}$$\end{definition}

\subsection{Proof of Theorem~\ref{thm-1.9}~(1,2)}
We suppose $f_{yy}>0$. The distributions
$$\ker(\mathfrak{R})=\operatorname{Span}\{\partial_{\tilde x}\}\text{ and }
\operatorname{Range}(\mathfrak{R})=\operatorname{Span}\{\partial_y,\partial_{\tilde x}\}
$$
are invariantly defined. To preserve these distributions, we set:
\begin{equation}\label{eqn-6.a}
\xi_1=a_{11}(\partial_x+f\partial_{\tilde x}+a_{12}\partial_y+a_{13}\partial_{\tilde x}),\quad
    \xi_2=\partial_y+a_{23}\partial_{\tilde x},\quad \xi_3=a_{33}\partial_{\tilde x}\,.
\end{equation}
To ensure that the inner products are normalized properly, we impose the relations:
$$a_{12}^2+2a_{13}=0,\quad a_{12}+a_{23}=0,\quad a_{11}a_{33}=1\,.$$
This determines $a_{13}$, $a_{23}$, and $a_{33}$; these parameters play no further role
and $\{\lambda,a_{11},a_{12}\}$ remain as free parameters where $\lambda$
is the \KV\ rescaling factor.
We suppose $f_{yyy}\ne0$. Set:
\begin{equation}\label{eqn-6.b}
\lambda:=f_{yyy}f_{yy}^{-1},\quad
a_{12}:=-f_{xyy}f_{yyy}^{-1},\quad
a_{11}^2:=f_{yy}^{-1}\lambda^2\,.
\end{equation}
We then have
\begin{eqnarray}
&&R(\xi_1,\xi_2,\xi_2,\xi_1)=a_{11}^2f_{yy}=\lambda^2,\nonumber\\
&&\nabla R(\xi_1,\xi_2,\xi_2,\xi_1;\xi_1)=a_{11}^3\{f_{xyy}+a_{12}f_{yyy}\}=0,\label{eqn-6.cx}\\
&&\nabla R(\xi_1,\xi_2,\xi_2,\xi_1;\xi_2)=a_{11}^2f_{yyy}=\lambda^2f_{yy}^{-1}f_{yyy}=\lambda^3\,.\nonumber
\end{eqnarray}
All the parameters of the theory have been determined (modulo a possible sign ambiguity in $a_{11}$) and
any \KV\ 1-model for $\mathcal{M}_f$ is isomorphic to $\mathfrak{M}_{1,0,1}$ in this special case. This proves Theorem~\ref{thm-1.9}~(1).

Suppose $f_{yy}>0$, $f_{yyy}=0$, and $f_{xyy}$ never vanishes. Set $f_{yy}=\alpha(x)$.
The parameter $a_{12}$ plays no role. To ensure that $\mathcal{M}_f$ is \KV\ 1-curvature
homogeneous, we impose the following relations where $\{a_{11},\lambda\}$ are unknown functions
to be determined and where $\{c_0,c_1\}$ are unknown constants:
\begin{eqnarray*}
&&R(\xi_1,\xi_2,\xi_2,\xi_1)=a_{11}^2(x)\alpha(x)=\lambda^2(x)c_0,\\
&&R(\xi_1,\xi_2,\xi_2,\xi_1;\xi_1)=a_{11}^3(x)\alpha_x(x)=\lambda^3(x)c_1,\\
&&R(\xi_1,\xi_2,\xi_2,\xi_1;\xi_2)=0\,.
\end{eqnarray*}
Consequently, $a_{11}^6(x)\alpha^3(x)=\lambda^6(x)c_0^3$ and
$a_{11}^6(x)\alpha^2_x(x)=\lambda^6(x)c_1^2$.
This shows that $\alpha^3(x)=c_3\alpha_x^2(x)$ for some constant $c_3$.
We solve this ordinary differential equation to see that
$$\alpha(x)=a(x-x_0)^{-2}\text{ for }0\ne a\in\mathbb{R}
\text{ and }x_0\in\mathbb{R}\,.$$
This has the form given in Theorem~\ref{thm-1.8}~(2b) and defines a locally homogeneous example. Theorem~\ref{thm-1.9}~(2) now follows.
\hfill\qed

\subsection{The proof of Theorem~\ref{thm-1.9}~(3)}
We assume that
$f_{yy}$ and $f_{yyy}$ never vanish as this case is the
only possible source of new examples not covered by
Theorem~\ref{thm-1.8}.
We shall suppose $f_{yy}>0$; the case $f_{yy}<0$ is handled similarly.
As any two \KV\ 1-curvature models for $\mathcal{M}_f$ are isomorphic,
we can adopt the normalizations of Equation~(\ref{eqn-6.a}),
(\ref{eqn-6.b}), and (\ref{eqn-6.cx}). We have:
\begin{eqnarray*}
&&R(\xi_1,\xi_2,\xi_2,\xi_1)=a_{11}^2f_{yy}=\lambda^2,\\
&&\nabla R(\xi_1,\xi_2,\xi_2,\xi_1;\xi_2)=a_{11}^2f_{yyy}=\lambda^3,\\
&&\nabla^2R(\xi_1,\xi_2,\xi_2,\xi_1;\xi_2,\xi_2)=a_{11}^2(x)f_{yyyy}=\lambda^4c_{11},\\
&&\frac{f_{yy}\cdot f_{yyyy}}{f_{yyy}\cdot f_{yyy}}=
\frac{\lambda^2a_{11}^{-2}\cdot\lambda^4c_{11}a_{11}^{-2}}{\lambda^6a_{11}^{-4}}=c_{11}\,.
\end{eqnarray*}
Thus $c_{11}$ is an invariant of the theory; this will imply the 3 families of the theory fall into
different local isometry types. The ordinary differential equation
$\frac{\alpha\alpha^{\prime\prime}}{\alpha^\prime\alpha^\prime}=c_{11}$
has the solutions $\alpha>0$ (see, for example, Lemma 1.5.5 of  \cite{G07}) of the form:
$$\alpha(t)=\left\{\begin{array}{ll}
e^{a(t+b)}&\text{ if }c_{11}=1\\
a(t+b)^c&\text{ if }c_{11}\ne1\end{array}\right\}
\text{ for }a\ne0\text{ and }c\ne0\,.$$
Thus
\begin{equation}\label{eqn-6.c}
f_{yy}=\left\{\begin{array}{ll}
e^{\alpha(x)(y+\beta(x))}&\text{ if }c_{11}=1\text{ for }\alpha(x)\ne0\\
\alpha(x)(y+\beta(x))^c&\text{ if }
c_{11}\ne1\text{ for }\alpha(x)\ne0\text{ and }c\ne0
\end{array}\right\}\,.\end{equation}

We wish to simplify Equation~(\ref{eqn-6.c}) to take $\beta(x)=0$.
We consider the change of variables
$T(x,y,z)=(x,y-\beta(x),\tilde x+y\beta_x(x))$:
$$\begin{array}{lrrr}
T_*\partial_x=(&1,&-\beta_x(x),&y\beta_{xx}(x)),\\
T_*\partial_y=(&0,&1,&\beta_x(x)),\\
T_*\partial_{\tilde x}=(&0,&0,&1).
\end{array}$$
We compute:
\begin{eqnarray*}
&&g_f(T_*\partial_x,T_*\partial_x)
=-2f(x,y-\beta(x))+\beta_x^2(x)+2y\beta_{xx}(x),\\
&&g_f(T_*\partial_x,T_*\partial_y)=0,
\qquad g_f(T_*\partial_x,T_*\partial_{\tilde x})=1,\\
&&g_f(T_*\partial_y,T_*\partial_y)=1,\qquad
g_f(T_*\partial_y,T_*\partial_{\tilde x})=g_f(T_*\partial_{\tilde x},T_*\partial_{\tilde x})=0.
\end{eqnarray*}
Thus $T^*g_f=g_{\tilde f}$ where
$$\tilde f(x,y)=f(x,y-\beta(x))-\textstyle\frac12
\left\{\beta_x^2(x)+2y\beta_{xx}(x)\right\}\,.$$
Consequently,
$\tilde f_{yy}(x,y)=f_{yy}(x,y-\beta(x))$.
Thus we may assume henceforth that $\beta(x)=0$ in
Equation~(\ref{eqn-6.c}), i.e.
$$
f_{yy}=\left\{\begin{array}{ll}
e^{\alpha(x)y}&\text{ if }c_{11}=1\text{ for }\alpha(x)\ne0\\
\alpha(x)y^c&\text{ if }c_{11}\ne1\text{ for }\alpha(x)\ne0\text{ and }c\ne0
\end{array}\right\}$$
We examine these two cases seriatum. We shall use the relations:
\begin{eqnarray}
\lambda&=&f_{yyy}f_{yy}^{-1},\quad
a_{12}=-f_{xyy}f_{yyy}^{-1},\quad
\lambda^2a_{11}^{-2}=f_{yy},\label{eqn-6.d}\\
\lambda^4c_{12}&=&\nabla^2R(\xi_1,\xi_2,\xi_2,\xi_1;\xi_1,\xi_2)=
a_{11}^3\{f_{xyyy}+a_{12}f_{yyyy}\},\label{eqn-6.e}\\
\lambda^4c_{11}&=&\nabla^2R(\xi_1,\xi_2,\xi_2,\xi_1;\xi_1,\xi_1)\label{eqn-6.f}\\
&=&a_{11}^4\{f_{xxyy}+2a_{12}f_{xyyy}+a_{12}^2f_{yyyy}-f_yf_{yyy}\}\,.
\nonumber\end{eqnarray}
\medbreak\noindent{\bf Case I.} Suppose $f_{yy}=e^{\alpha(x)y}$.
Then Equation~(\ref{eqn-6.d}) implies:
$$
\lambda=f_{yyy}f_{yy}^{-1}=\alpha(x),\
 a_{12}=-f_{xyy}f_{yyy}^{-1}=-\alpha_x(x)\alpha(x)^{-1},\
\lambda^2a_{11}^{-2}=f_{yy}=e^{\alpha(x)y}\,.
$$
We use Equation~(\ref{eqn-6.e})
to see that:
\begin{eqnarray*}
&&f_{xyyy}+a_{12}f_{yyyy}=\partial_x\{\alpha(x)e^{\alpha(x)y}\}
-\alpha_x(x)\alpha(x)e^{\alpha(x)y}\\
&=&
\alpha_x(x)\cdot e^{\alpha(x)y}=a_{11}^{-3}\lambda^4c_{12}=\alpha(x)e^{\frac32\alpha(x)y}c_{12}\,.
\end{eqnarray*}
It now follows that $\alpha_x(x)=0$ so $\alpha(x)=a$ is constant. Thus
we may express:
$$f(x,y)=a^{-2}e^{ay}+u(x)y+v(x)\,.$$
We then use Equation~(\ref{eqn-6.d}) to see
$$\lambda=a,\quad a_{12}=0,\quad \lambda^2a_{11}^{-2}=e^{ay}\,.$$
Equation~(\ref{eqn-6.f}) then leads to the identity:
$$e^{2ay}c_{11}=a_{11}^{-4}\lambda^4c_{11}=-f_yf_{yyy}=-e^{2ay}-u(x)ae^{ay}\,.$$
This implies that $u(x)=0$ and hence $f=a^{-2}e^{ay}+v(x)$. Let $w_x=v(x)$ and set:
$$\begin{array}{rrrr}
T(x,y,\tilde x)=(&x,&y,&\tilde x+2w(x)),\\
T_*\partial_x=(&1,&0,&2v(x)),\\
T_*\partial_y=(&0,&1,&0),\\
T_*\partial_{\tilde x}=(&0,&0,&1).
\end{array}$$
Under this change of variables:
\begin{eqnarray*}
&&g(T_*\partial_x,T_*\partial_x)
=-2a^{-2}e^{ay}-2v(x)+2v(x)=-2a^{-2}e^{ay},\\
&&g(T_*\partial_x,T_*\partial_y)=g(T_*\partial_y,T_*\partial_{\tilde x})
=g(T_*\partial_{\tilde x},T_*\partial_{\tilde x})=0,\\
&&g(T_*\partial_x,T_*\partial_{\tilde x})=g(T_*\partial_y,T_*\partial_y)=1\,.
\end{eqnarray*}
Thus we may take $f=a^{-2}e^{ay}$. Replacing $y$ by $y+y_0$ for
suitably chosen $y_0$, then replaces $f$ by $e^{ay}$ as desired.

We now show $\mathcal{M}_{e^{ay}}$ is a homogeneous space. Set:
$$T(x,y,\tilde x)=(\pm e^{-ay_0/2}x+x_0,y+y_0,\pm e^{ay_0/2}\tilde x
+\tilde x_0)\,.$$
Then
$$T_*\partial_x=\pm e^{-ay_0/2}\partial_x,\quad
T_*\partial_y=\partial_y,\quad
T_*\partial_{\tilde x}=\mp e^{ay_0/2}\partial_{\tilde x}\,.
$$
We show
that $T$ is an isometry by verifying:
\begin{eqnarray*}
&&g(T_*\partial_x,T_*\partial_x)(y+y_0)=-2e^{-ay_0}e^{a(y+y_0)}
=g(\partial_x,\partial_x)(y),\\
&&g(T_*\partial_x,T_*\partial_y)=0,\quad
g(T_*\partial_x,T_*\partial_{\tilde x})=1,\quad
g(T_*\partial_y,T_*\partial_y)=1,\\
&&g(T_*\partial_y,T_*\partial_{\tilde x})=0,\quad
g(T_*\partial_{\tilde x},T_*\partial_{\tilde x})=0.
\end{eqnarray*}
Since $(x_0,y_0,\tilde x_0)$ are arbitrary,
$\mathcal{I}(\mathcal{M}_{e^{ax}})$
acts transitively on $\mathbb{R}^3$ so this manifold is globally homogeneous.
This verifies Theorem~\ref{thm-1.9}~(3a).

\medbreak\noindent{\bf Case II.}
Suppose that $f_{yy}=\alpha(x)y^c$ for $\alpha(x)>0$ and $c\ne0$.
Equation~(\ref{eqn-6.d}) yields
$$
\lambda=f_{yyy}f_{yy}^{-1}=cy^{-1},\quad
a_{12}=-f_{xyy}f_{yyy}^{-1}=-\textstyle\frac{a_x(x)y}{c\alpha(x)},
\quad\lambda^2a_{11}^{-2}=f_{yy}=\alpha(x)y^c\,.
$$
We apply Equation~(\ref{eqn-6.e}) to see:
\begin{eqnarray*}
&&\textstyle f_{xyyy}+a_{12}f_{yyyy}=\alpha_x(x)cy^{c-1}
-\frac{\alpha_x(x)y}{c\alpha(x)}c(c-1)\alpha(x)y^{c-2}
=\alpha_x(x)y^{c-1}\\
&=&a_{11}^{-3}\lambda^4c_{12}=\alpha(x)^{3/2}y^{3/2c}cy^{-1}c_{12}\,.
\end{eqnarray*}
This implies
$$\alpha_x(x)\alpha(x)^{-3/2}=c\cdot c_{12}\cdot y^{c/2}\,.$$
Consequently $\alpha_x(x)=0$ so $\alpha(x)=a$ is constant. Consequently,
$f_{yy}=ay^c$ for $c\ne0$ and $a\ne0$.
Let $P(t)$ solve the equation $P^{\prime\prime}(t)=t^c$. We then have
$$f(y)=aP(y)+u(x)y+v(x)\,.$$
We apply Equation~(\ref{eqn-6.f}) with $a_{12}=0$:
\begin{eqnarray*}
&&a_{11}^{-4}\nabla^2R(\xi_1,\xi_2,\xi_2,\xi_1;\xi_1,\xi_1)=
-f_yf_{yyy}=-aP^\prime(y)acy^{c-1}-u(x)acy^{c-1}\\
&=&c_{11}\lambda^4a_{11}^{-4}=c_{11}a^2y^{2c}\,.
\end{eqnarray*}
If $c=-1$, then $P^\prime(y)=\ln(y)$ and this relation is impossible. Consequently $c\ne-1$ and
we may conclude that $u(x)=0$.
We therefore have $f=aP(y)+v(x)$. As in Case I, the constant term is
eliminated and $a$ is set to $1$ by making a change
of variables
$$T(x,y,\tilde x)=(a^{-1/2}x,y,a^{1/2}\tilde x+2w(x))$$
where $w_x(x)=v(x)$. Thus $f=\pm\ln(y)$ or $f=\pm y^\varepsilon$ for
$\varepsilon\ne0,1,2$.

\medbreak\noindent{\bf Case II-a.} Let $f(y)=\ln(y)$; the case
$f(y)=-\ln(y)$ is similar. We know by Theorem~\ref{thm-1.8}
that $\mathcal{M}_f$ is not $2$-curvature homogeneous and hence
is not homogeneous. For $\lambda>0$ and $(x_0,\tilde x_0)$ arbitrary, set:
$$T(x,y,\tilde x):=(\lambda x+x_0,\lambda y,\lambda\tilde x
+\tilde x_0+\lambda\ln\lambda x)\,.$$
We compute:
\begin{eqnarray*}
&&T_*\partial_x=\lambda\partial_x+\lambda\ln\lambda\partial_{\tilde x},
\quad T_*\partial_y=\lambda\partial_y,\quad
T_*\partial_{\tilde x}=\lambda\partial_{\tilde x},\\
&&g(T_*\partial_x,T_*\partial_x)(\lambda y)
=\lambda^2\{-2\ln(y)-2\ln\lambda\}+2\lambda^2\ln\lambda
=\lambda^2g(\partial_x,\partial_x)(y),\\
&&g(T_*\partial_x,T_*\partial_y)=0,\quad
g(T_*\partial_x,T_*\partial_{\tilde x})=\lambda^2,\quad
g(T_*\partial_y,T_*\partial_y)=\lambda^2,\\
&&g(T_*\partial_y,T_*\partial_{\tilde x})=0,\quad
g(T_*\partial_{\tilde x},T_*\partial_{\tilde x})=0\,.
\end{eqnarray*}
This defines a transitive action on
$\mathbb{R}\times\mathbb{R}^+\times\mathbb{R}$.

\medbreak\noindent{\bf Case II-b.} Let $f(y)=y^c$ for $c\ne0,1,2$.
Again, Theorem~\ref{thm-1.8} implies $\mathcal{M}_f$ is not
$2$-curvature homogeneous and hence not homogeneous. Let
$$T(x,y,\tilde x)=(\lambda^{(2-c)/2}x+x_0,\lambda y,
\mp\lambda^{2+(c-2)/2}\tilde x+\tilde x_0)\,.$$
We compute:
\begin{eqnarray*}
&&T_*\partial_x=\lambda^{(2-c)/2}\partial_x,\quad
T_*\partial_y=\lambda\partial_y,\quad
T_*\partial_{\tilde x}=\lambda^{2+(c-2)/2}\partial_{\tilde x},\\
&&g(T_*\partial_x,T_*\partial_x)(\lambda y)=
-2\lambda^{(2-c)}\lambda^cy^c=\lambda^2g(\partial_x,\partial_x)(y),\\
&&g(T_*\partial_x,T_*\partial_y)=0,\quad
g(T_*\partial_x,T_*\partial_{\tilde x})=\lambda^2,\quad
g(T_*\partial_y,T_*\partial_y)=\lambda^2,\\
&&g(T_*\partial_y,T_*\partial_{\tilde x})=0,\quad
g(T_*\partial_{\tilde x},T_*\partial_{\tilde x)}=0\,.
\end{eqnarray*}
Thus $T$ is a homothety; since $\lambda>0$ is arbitrary and since
$(x_0,\tilde x_0)$ are arbitrary, the group of homotheties acts transitively
on $M$. This verifies Theorem~\ref{thm-1.9}~(3c) and completes
the proof of Theorem~\ref{thm-1.9}.\hfill\qed

\section{The proof of Theorem~\ref{thm-1.11}}\label{sect-7}

We consider the family of Lorentzian manifolds
given in Section~\ref{sect-6} and consider the family of models where
$\langle\xi_1,\xi_3\rangle=\langle\xi_2,\xi_2\rangle=1$
defines the inner product on the vector space
$V=\operatorname{Span}\{\xi_1,\xi_2,\xi_3\}$.
Choose the isometry between $T_PM$ to be defined by:
$$\xi_1=a_{11}(\partial_x+f\partial_{\tilde x}+a_{12}\partial_y+a_{13}\partial_{\tilde x}),\quad
\xi_2=\partial_y+a_{23}\partial_{\tilde x},\quad
\xi_3=a_{11}^{-1}\partial_{\tilde x}\,.
$$
To ensure the map in question is an isometry, we require
$a_{12}^2+2a_{13}=0$ and $a_{12}+a_{23}=0$.
This normalizes the parameters $a_{13}$ and $a_{23}$; $a_{11}$ and $a_{12}$ and $\lambda$ are the free parameters
of the theory.

We take $f(x,y)=\frac12\alpha(x)y^2$. The only (possibly)
non-zero covariant derivatives are given by:
$$
\nabla^\ell R(\partial_x,\partial_y,\partial_y,\partial_x;\partial_x,...,\partial_x)=\alpha^{(\ell)}(x)\,.
$$
The parameter $a_{12}$ plays no role in the theory and we have
$$\nabla^\ell R(\xi_1,\xi_2,\xi_2,\xi_1;\xi_1,...,\xi_1)=a_{11}^{2+\ell}\alpha^{(\ell)}\,.$$

Suppose $k$ is given. Let $\alpha^{(k)}:=e^{-x^2}$ and recursively define
\begin{equation}\label{eqn-7.a}
\alpha^{(\ell)}(x):=\int_{t=-\infty}^x\alpha^{(\ell+1)}(t)dt\text{ for }0\le\ell\le k-1\,.
\end{equation}
For $x\le -1$, $-x^2\le x$ and consequently $\alpha^{(k)}\le e^{x}$.
Inductively, suppose $\alpha^{(\ell)}\le e^{x}$
for $x\le-1$. Then integration yields $\alpha^{(\ell-1)}\le e^x$ for $x\le-1$ as well. Thus the
integrals in Equation~(\ref{eqn-7.a}) converge and $\alpha$ is a smooth function.
We have $\alpha^{(\ell)}>0$ for $0\le\ell\le k$; setting $a_{11}^{\ell}=\{\alpha^{(\ell)}\}^{1/(2+\ell)}$ then yields
$$\nabla^\ell R(\xi_1,\xi_2,\xi_2,\xi_1;\xi_1,...,\xi_1)=1$$
and shows $\mathcal{M}_\alpha$ is variable $k$-curvature homogeneous. Since, however,
$\alpha^{(k+1)}$ vanishes when $x=0$ and is non-zero when $x\ne0$,
$\mathcal{M}_\alpha$ is
not variable $k+1$-curvature homogeneous on any neighborhood
of $0$. This establishes
Theorem~\ref{thm-1.11}~(1).

We take $\alpha(x)=e^x$. Then $\alpha^{(k)}(x)>0$ for all $x$ and for all $k$. Thus
the argument given above shows $\mathcal{M}_\alpha$ is variable $k$-curvature homogeneous
for all $k$. On the other hand, by Theorem~\ref{thm-1.9}~(2),
$\mathcal{M}_\alpha$ is not \KV\ 1-curvature homogeneous.
This establishes Theorem~\ref{thm-1.11}~(2).

We consider cases to establish Theorem~\ref{thm-1.11}~(3):
\begin{enumerate}
\item We use Theorem~\ref{thm-1.8}~(1b)
and Theorem~\ref{thm-1.9}~(1) to see that
the implication $(1b)\Rightarrow(1a)$ fails. This also shows the implication $(2b)\Rightarrow(1a)$ fails.
\item If we take $f=\frac12e^xy^2$, then
this is variable $k$-curvature homogeneous for all $k$ by Theorem~\ref{thm-1.11}~(2)
but not, by Theorem~\ref{thm-1.9}~(2) \KV\ 1-curvature homogeneous.
Thus the implications
$(2a)\Rightarrow(1a)$, $(2a)\Rightarrow(1b)$, $(2b)\Rightarrow(1a)$, and $(2b)\Rightarrow(1b)$
fail.
\item Let $(S^2,g_2)$ denote the standard round sphere in $\mathbb{R}^3$. Let
$$\mathcal{M}_t:=(M,g_t):=(\mathbb{R}\times S^2,e^{tx}(dx^2+g_2))$$
where $t\in(0,\epsilon)$ is a small positive real
parameter. We showed previously that this is homothety homogeneous
and hence satisfies (1b) for all $k$. We also showed it was not
$0$-curvature homogeneous for generic values of $t$. Note that
variable $0$-curvature and $0$-curvature homogeneous are the same.
Thus the implications $(1b)\Rightarrow(2a)$ and $(2b)\Rightarrow(2a)$ fail.
\end{enumerate}

\section*{Acknowledgments}
Research of all the authors was partially supported by project MTM2009-07756. Research of P. Gilkey and S. Nik\v cevi\'c was also
partially supported by project 174012 (Serbia).

\end{document}